\documentstyle{article}
\input epsf

\newtheorem{theorem}{Theorem}

\def\G{\Gamma}

\def\la{\lambda}

\title{$S$--Partitions}
\date{}
\author{William M. Y. Goh
\\
{\small  Department of Mathematics and Computer Science}\\
{\small Drexel University}\\
{\small Philadelphia, PA 19104}\\
{\small \texttt{wgoh@mcs.drexel.edu}}\\
\and Pawe{\l} Hitczenko \\
{\small  Department of Mathematics and Computer Science}\\
{\small Drexel University}\\
{\small Philadelphia, PA 19104}\\
{\small \texttt{phitczen@mcs.drexel.edu}}\\
\and Ali Shokoufandeh \\
{\small  Department of Mathematics and Computer Science}\\
{\small Drexel University}\\
{\small Philadelphia, PA 19104}\\
{\small \texttt{ashokouf@mcs.drexel.edu}}
}

\begin{document}

\begin{titlepage}

\maketitle

\abstract{This  note reports  on  the number  of  $s$-partitions of  a
  natural number  $n$. In  an $s$-partition of  $n$ each cell  has the
  form $2^k-1$  for some integer $k$.  Such  partitions have potential
  applications   in    cryptography,   specifically   in   distributed
  computations  of  the  form  $a^n  {\  \rm  mod\  }  m$.   The  main
  contribution of  this paper  is a correction  to the  upper-bound on
  number  of $s$-partitions  presented by  Bhatt~\cite{bha}.   We will
  give a precise  asymptotics for the number of  such partitions for a
  given integer $n$.

\vspace{1cm}
\noindent {\em Keywords}: Combinatorial problems, integer partitions}
\end{titlepage}
\section{Introduction}
The $s$-partition of an integer  $n$ is a decomposition $n=\sum_i n_i$
such that  each $n_i$ is  of the form  $2^k-1$, for some integer  $k$.
Closely related  are {\em  binary partitions}, i.e.   partitions whose
cells are powers of 2 (see~\cite{pen}).  Perhaps, the simplest form of
a binary partition of a  number is its binary representation, in which
each part  size has the  form $2^k$, for  some integer $k$.   Given an
integer  $n$  its  binary  representation  can be  computed  by  using
consecutive shifts in $\Theta(\log n)$ steps.

An important computation  in most of the cryptographic  systems is the
computations of  the form  $a^n {\  \rm mod\ }  m$ (see~\cite{sch}).
In~\cite{bha},  Bhatt presented  an NC-algorithm  for  computations of
this form on  a CRCW parallel model.  While  giving a characterization
of  $s$-partitions, for  a given  integer $n$,  he presented  a simple
algorithm for  computing its $d$-partition
in  $\Theta(\log  n)$.   He  went  on to  use  this  decomposition  in
computing  $a^n  {\  \rm mod\  }  m$  in  poly-logarithmic time  on  a
distributed system, using polynomial  number of processors in terms of
$n$.

As  the main theoretical
result Bhatt gave  the quantity
$$2+\left \lfloor {n\over 3}\right \rfloor+ \sum_{i=0}^{\lfloor \log n
  \rfloor}{\lfloor     \log(n-3i)\rfloor}^     {\lfloor     \log(n-3i)
  \rfloor-1},$$
(where $\log$ means $\log_2$) as an upper bound for the number of $s$ partitions of an integer $n$,
but  this bound is not correct.
The purpose of  this note is to present a correction  to his result by
giving a  precise asymptotics for  the number of
$s$-partitions. Specifically, we have
\begin{theorem}
Let $p_s(n)$ denote the total  number of $s$--partitions of an integer
$n$     and     for     a     function     $f(x)=\left\lfloor\frac{\ln
x}{\ln2}\right\rfloor-\frac{\ln      x}{\ln2}+\frac12$,      let     $
\alpha:=\lim_{u\to\infty}\int_2^u\frac{f(v)}{v(v-1)}dv=-0.4934\dots.\label{alpha}$
Then,
\begin{eqnarray*}
\ln p_s(n)&=&\frac1{2\ln2}(\ln (n+1)-\ln\ln (n+1)+\ln
 \ln2)^2 \\
& &\qquad\qquad+(\frac1{\ln2}-\frac32)\ln (n+1)+\ln\ln (n+1)-\ln\ln2 \\
&+&W(\ln (n+1)-\ln\ln (n+1)-\ln \ln2)-\frac12\ln 2\pi+H+o(1),
\end{eqnarray*}
where
$$H=\frac{\pi^2+\ln^22}{12\ln2}+\alpha
+\frac1{\ln2}\int_0^\infty\frac{\ln v-\ln(1-e^{-v})}{e^v-1}dv,$$
and
$$
W(z)=-\sum_{\nu\ne0}\left(\frac{2\pi\nu}{\ln2}\right)^2\G\left(\frac{2\pi
    i\nu}{\ln2}\right)\zeta\left(1+\frac{2\pi   i\nu}{\ln2}\right)c_\nu
e^{(2\pi i \nu/\ln2)z}.
$$
\end{theorem}
A proof is a quick consequence  of   a   theorem
that
was proven with binary partitions in mind. Such  partitions have been studied
by several
authors beginning in the  40's. Not surprisingly, the results obtained
by these authors are general enough to be applicable in the context of
$s$-partitions.   We   will   recall   one such   theorem proved in
\cite{pen}     in the next section.
\section{Proof of Theorem~1}
\label{s-par}
Here is a relevant  part of a  theorem proved by
Pennington. We refer the reader  to his paper~\cite{pen} for more
details and credits concerning pre-1953 activities.
\begin{theorem} Let $0<\la_1<\la_2<\dots$ be a given sequence of numbers
  with
  $$N(u)=\sum_{\la_\nu\le u}1=a\ln u+b+R(u)$$
  for $u>0$, where
  $$\int_{\la_1}^u\frac{R(v)}{v}dv=c+V\left(\frac{\ln
      u}{\rho}\right)+o(1)$$
  as  $u\to\infty$,  $a$,  $b$,  $c$,  $\rho$  are  constants,  $a>0$,
  $\rho>0$, and $V$ is a  periodic function with period 1, bounded and
  integrable   in  the   interval  $0\le   x\le  1$.   Let  $\{c_\nu:\
  -\infty<\nu<\infty\}$ be the complex Fourier coefficients of $V$ and
  suppose  that $c_0=0$  and  $\sum_{\nu\ne0}|c_\nu/\nu|<\infty$.  For
  real $u$ let $P(u)$ be the number of solutions of the inequality
  $$r_1\la_1+r_2\la_2+r_3\la_3+\dots<u$$
  in  integers $r_\nu\ge0$, and
  let $P_h(u)=\{P(u)-P(u-h)\}/h$.  Then, if $h$ is a positive constant
  for which $P_h(u)$ is an  increasing function of $u$ (this condition
  is certainly satisfied if $h$ belongs to the sequence $\la_\nu$), as
  $u\to\infty$,
\begin{eqnarray*}
\ln P_h(u)&=&\frac12a(\ln u-\ln\ln u-\ln a)^2+(a-\frac12)\ln u \\
& &\qquad\qquad+(b-\frac12)(\ln u-\ln\ln u-\ln a) \\
&+&W(\ln u-\ln\ln u-\ln a)-\frac12\ln 2\pi+H+o(1),
\end{eqnarray*}
where
$$H=c-b\ln\la_1-\frac12a\ln^2\la_1+a\int_0^\infty\frac{\ln
  v-\ln(1-e^{-v})}{e^v-1}dv,$$
and
$$
W(z)=-\sum_{\nu\ne0}\left(\frac{2\pi\nu}{\rho}\right)^2\G\left(\frac{2\pi
    i\nu}{\rho}\right)\zeta\left(1+\frac{2\pi   i\nu}{\rho}\right)c_\nu
e^{(2\pi i \nu/\rho)z}.
$$
The function  $W(z)$  is defined  for  all $z=x+iy$  in the  strip
$|y|<\frac12\pi$, and is bounded and uniformly continuous in any fixed
interior  strip  $|y|\le  M<\frac12\pi$.
\end{theorem}
We  observe that  if $\la_1=1$  then, for  an integer  $n$, $P_1(n+1)$
counts the number of partitions of  $n$ into parts of sizes in the set
$\{\la_\nu:\  \nu\ge1\}$.  In  particular,  $p_s(n)=P_1(n+1)$ for  the
sequence  $\la_\nu=2^\nu-1$,  $\nu\ge1$.   Thus,  in  order  to  prove
Theorem~1,  it suffices  to  check that  this  sequence satisfies  the
assumptions of Pennington's theorem with
\begin{equation}
a=\frac1{\ln2}, \quad b=-\frac12,\quad
c=\frac{\pi^2+\ln^22}{12\ln2}+\alpha,\quad \rho=\ln2,\quad
\label{constants}
\end{equation}
   To this end we write
\begin{equation}
N(u)           =          \sum_{2^\nu\le          u+1\atop
  \nu\ge1}\!\!\!\!\!1=\left\lfloor\frac{\ln(u+1)}{\ln2}\right\rfloor
:= \frac{\ln u}{\ln2}-\frac12+R(u), \label{eq1}
\end{equation}
and we will show that
$$R(u):=\frac{\ln(1+1/u)}{\ln2}+
\left(\left\lfloor\frac{\ln(u+1)}{\ln2}\right\rfloor-\frac{\ln(u+1)}{\ln2}+\frac12\right)
$$
has the desired properties. Since $\la_1=1$, we have
\begin{equation}
\int_{\la_1}^u\frac{R(v)}{v}dv=\int_1^u\frac{\ln(1+1/v)}{\ln2}\frac{dv}v+\int_
1^u\frac{f(v+1)}{v}dv,\label{eq2}
\end{equation}
where   $f(x)=\left\lfloor\frac{\ln
    x}{\ln2}\right\rfloor-\frac{\ln  x}{\ln2}+\frac12$ as defined in
Theorem~1.   Changing the
variables $t=1/v$ in the first integral we see that it converges to
\begin{equation}
\frac1{\ln2}\int_0^1\frac{\ln(1+t)}{t}dt=\frac{\pi^2}{12\ln2}.\label{eq3}
\end{equation}
To handle the second integral, we rely on work of Pennington who showed
that
\begin{equation}
\int_1^u\frac{f(v)}vdv=\frac{\ln2}{12}-\sum_{\nu\ne0}\frac{\ln2}{4\pi^2\nu^2}e^{(2\pi
i\nu/\ln2)\ln u} .\label{eq4}
\end{equation}
The difference between the two is
\begin{eqnarray*}
& &\int_1^u\frac{f(v+1)}{v}dv-\int_1^u\frac{f(v)}{v}dv=\int_2^{u+1}
\frac{f(v)}{v-1}dv-\int_1^u\frac{f(v)}{v}dv \\
&=&
\int_2^u\frac{f(v)}{v(v-1)}dv-\int_1^2\frac{f(v)}vdv+\int_u^{u+1}\frac{f(v)}{v-1}dv.
\end{eqnarray*}
An elementary calculation  shows that the middle integral  is zero and
since  the  function  $f$  is  uniformly bounded,  the  last  integral
vanishes as $u\to \infty$.  The first integral converges absolutely to
$\alpha$   defined  in   Theorem~1.   Comparing   (\ref{eq1})  through
(\ref{eq4})  with  the  statement   of  Theorem~2,  we  see  that  its
conditions are satisfied with the  constants $a$, $b$, $c$, and $\rho$
given by (\ref{constants}) and  with Fourier coefficients $ c_0=0$ and
$c_\nu=-\frac{\ln2}{4\pi^2\nu^2}$    for   $\nu\ne0$.    This   proves
Theorem~1.


\begin{thebibliography}{WWW}


\bibitem{bha} Bhatt,  P.C.P., {\em An  interesting way to  partition a
    number},  Inform.  Process.  Lett.,  {\bf 71}  (1999),  141--148.
  (http://sciencedirect.com)


\bibitem{pen}  Pennington, W.B.,  {\em On  Mahler  partition problem},
  Ann.Mathematics {\bf 57} (1953), 579--589 (http://links.jstor.org).

\bibitem{sch} Schneier, B., {\em Applied Cyptography}, John Wiley, New York,
  1996.

\end{thebibliography}
\end{document}